# $\gamma_0$-compact , $\gamma^s$-regular and $\gamma^s$-normal spaces

**Bashir Ahmad*** **and Sabir Hussain** **

* Centre for Advanced studies in Pure and Applied Mathematics, Bahauddin Zakariya University, Multan, Pakistan.
Present Address: Department of Mathematics, King Abdul Aziz University, P.O.Box 80203, Jeddah 21589, Saudi Arabia.
E.mail: drbashir9@gmail.com

** Department of Mathematics, Islamia University, Bahawalpur, Pakistan.
E.mail: sabiriub@yahoo.com

**Abstract**. We define and study the properties of $\gamma^*$-regular and $\gamma$-normal spaces. We also continue studying $\gamma_0$-compact spaces defined in [5].

**Introduction**. In 1979, S. Kasahara [7] defined an operation $\alpha$ on topological spaces. He introduced and studied $\alpha$-closed graphs of a function. In 1983, D.S. Jankovic [6] defined $\alpha$- closed set and further worked on functions with $\alpha$- closed graphs. In 1991, H. Ogata [8] introduced the notions of $\gamma$-$T_i$ ,i = 0,1/2,1,2; and studied some topological properties. In 1992 ( respt. in 1993) F.U.Rehman and B.Ahmad [9] ( [2] ) defined and investigated several properties of $\gamma$-interior, $\gamma$-exterior, $\gamma$-closure and $\gamma$-boundary points in topological spaces ( respt. in product spaces [2] ), and studied the characterizations of $(\gamma,\beta)$-continuous mappings initiated by H. Ogata [8]. In 2003, B.Ahmad and S.Hussain [3] continued studying the properties of $\gamma$-operations on topological spaces introduced by S. Kasahara [7]. Recently, B.Ahmad and S.Hussain [4] defined $\gamma$-nbd, $\gamma$-nbd base at x, $\gamma$-closed nbd, $\gamma$-limit point, $\gamma$-isolated point, $\gamma$-convergent point and $\gamma^*$-regular space and discussed their several properties. They further established the properties of $(\gamma,\beta)$- continuous, $(\gamma,\beta)$- open functions and $\gamma$-$T_2$ spaces.

In this paper, we continue studying $\gamma_0$-compact space defined in [5] and study its properties. We also define and study some properties of $\gamma^s$-regular and $\gamma^s$-normal spaces.

First, we recall some definitions and results used in this paper. Hereafter we shall write space in place of topological space.

**Definition** [8]. Let $(X,\tau)$ be a space. An operation $\gamma : \tau \to P(X)$ is a function from $\tau$ to the power set of X such that $V \subseteq V^\gamma$, for each $V \in \tau$, where $V^\gamma$ denotes the value of $\gamma$ at V. The operations defined by $\gamma(G) = G$, $\gamma(G) = cl(G)$ and $\gamma(G) = intcl(G)$ are examples of operation $\gamma$.

**Definition** [9]. Let $A \subseteq X$. A point $a \in A$ is said to be $\gamma$-interior point of A iff there exists an open nbd N of a such that $N^\gamma \subseteq A$ and we denote the set of all such points by $int_\gamma(A)$. Thus

$$int_\gamma(A) = \{ x \in A : x \in N \in \tau \text{ and } N^\gamma \subseteq A \} \subseteq A.$$

Note that A is $\gamma$-open [7] iff $A = int_\gamma(A)$.

A set A is called $\gamma$- closed [7] iff $X-A$ is $\gamma$-open.

**Definition** [9]. A point $x \in X$ is called a $\gamma$-closure point of $A \subseteq X$, if $U^\gamma \cap A \neq \emptyset$, for each open nbd U of x. The set of all $\gamma$-closure points of A is called $\gamma$-closure of A and is denoted by $cl_\gamma(A)$. A subset A of X is called $\gamma$-closed, if $cl_\gamma(A) \subseteq A$. Note that $cl_\gamma(A)$ is contained in every $\gamma$-closed superset of A.

**Definition 1**. An operation $\gamma : \tau \to P(X)$ is said to be strictly regular, if for any open nbds U , V of $x \in X$, there exists an open nbd W of x such that $U^\gamma \cap V^\gamma = W^\gamma$.

**Definition 2**. An operation $\gamma : \tau \to P(X)$ is said to be $\gamma$-open, if $V^{\gamma}$ is $\gamma$-open for each $V \in \tau$.

**Example 1**. Let $X=\{a,b,c\}$, $\tau =\{ \emptyset, X, \{a\}, \{b\},\{a,b\}\}$.

Define an operation $\gamma : \tau \to P(X)$ by $\gamma(A) = \text{intcl}(A)$.

Clearly the $\gamma$-open sets are only $\emptyset$, X, $\{a\}$, $\{b\}$,$\{a,b\}$. It is easy to see that $\gamma$ is strictly regular and $\gamma$-open on X.

**Example 2**. Let $X=\{a,b,c\}$, $\tau =\{ \emptyset, X, \{a\}, \{b\},\{a,b\}\}$.

Define an operation $\gamma : \tau \to P(X)$ by $\gamma(A) = cl(A)$.

Clearly the $\gamma$-open sets are only $\emptyset$, X. It is easy to see that $\gamma$ is strictly regular but not $\gamma$-open on X.

## 1. $\gamma_0$-compact space

**Definition** [5]. A space X is said to be $\gamma_0$-compact, if for every cover $\{V_i : i \in I\}$ of X by $\gamma$-open sets of X, there exists a finite subset $I_0$ of I such that $X= \bigcup_{i\in I_0} cl_{\gamma} ( V_i)$.

Then the following characterization of a $\gamma_0$-compact space is immediate :

**Theorem** [5] . A space X is $\gamma_0$-compact iff every class of $\gamma$-open and $\gamma$- closed sets with empty intersection has a finite subclass with empty intersection.

**Definition** [8]. A space X is said to be $\gamma$-$T_2$ space. If for each pair of distinct points x, y in X, there exist open sets U , V such that $x \in U$, $y \in V$ and $U^{\gamma} \cap V^{\gamma} = \emptyset$.

**Thorem 1**. Let X be a $\gamma$-$T_2$ space and suppose that C be a $\gamma_0$-compact subset of X and $x \in X - C$, then there are open sets $U_x$ and $V_x$ in X such that $x \in U^{\gamma}_x$ and $C \subseteq V^{\gamma}_x$ and $U^{\gamma}_x \cap V^{\gamma}_x = \emptyset$, where $\gamma$ is regular and $\gamma$-open.

**Proof**. Suppose that C is $\gamma_0$-compact subset of X and $x \in X - C$. For each $y \in C$ , $y \neq x$.

Since X is $\gamma$-$T_2$ , there are open sets $U_{xy}$ and $V_y$ containing x and y respectively such that $U^{\gamma}_{xy} \cap V^{\gamma}_{y} = \varnothing$. Now, let $\{ V^{\gamma}_{y} \cap C : y \in C \}$ be $\gamma$-open cover of C. Since C is $\gamma_0$-compact , then $\gamma$-open cover has a finite subset $\{ V^{\gamma}y_1 \cap C , V^{\gamma}y_2 \cap C, \ldots , V^{\gamma}y_n \cap C \}$ such that $C = \bigcup_{i=1}^{n} cl_{\gamma} ( V^{\gamma}y_i \cap C )$. Let $U^{\gamma}y_1 , U^{\gamma}y_2 , \ldots , U^{\gamma} y_n$ be the corresponding $\gamma$-open ses containing x. Take

$$U^{\gamma}_{x} = ( \bigcap_{i=1}^{n} cl_{\gamma} (U^{\gamma}xy_i))$$

and

$$V^{\gamma}_{x} = ( \bigcup_{i=1}^{n} cl_{\gamma} (V^{\gamma}xy_i)),$$

then $x \in U^{\gamma}_{x}$ and $C \subseteq V^{\gamma}_{x}$ . Where $U^{\gamma}_{x}$ and $V^{\gamma}_{x}$ are $\gamma$-closed , since $\gamma$ is regular.

Also

$$U^{\gamma}_{x} \cap V^{\gamma}_{x} = ( \bigcap_{i=1}^{n} cl_{\gamma} (U^{\gamma}xy_i )) \cap ( \bigcup_{i=1}^{n} cl_{\gamma} (V^{\gamma}y_i ))$$

$$= \bigcap_{i=1}^{n} ( \bigcup_{i=1}^{n} cl_{\gamma} (U^{\gamma}xy_i) \cap cl_{\gamma} (V^{\gamma}y_i ))$$

$$= \bigcap_{i=1}^{n} ( \bigcup_{i=1}^{n} cl_{\gamma} (U^{\gamma}xy_i \cap V^{\gamma}y_i )) \qquad ( \gamma \text{ is regular [9] } )$$

$$= \bigcap_{i=1}^{n} ( \bigcup_{i=1}^{n} cl_{\gamma} ( \varnothing )) = \varnothing.$$

**Theorem 2**. Let X be a $\gamma$-$T_2$ space. Then every $\gamma_0$-compact subset A of X is $\gamma$-closed, where $\gamma$ is regular and $\gamma$-open.

**Proof**. Let X be a $\gamma$-$T_2$ space and A be a $\gamma_0$-compact subset of X. We show that $X - A$ is $\gamma$-open. For this, let $x \in X - A$. Then $y \in A$ gives $x \neq y$. Since X is $\gamma$-$T_2$ , there are open sets $U_{xy}$ and $U_y$ in X containing x and y respectively such that $U^{\gamma}_{xy} \cap U^{\gamma}_{y} = \varnothing$.

Let the collection $\{ U^{\gamma}_{y} \cap A : y \in A \}$ is a cover of A by $\gamma$-open and $\gamma$-closed sets of A.

Since A is $\gamma_0$-compact, there is a finite subset $\{ U^{\gamma}y_i \cap A : i = 1,2, \ldots ,n\}$ such that

$$A = \bigcup_{i=1}^{n} cl_{\gamma} (U^{\gamma}y_i \cap A) = \bigcup_{i=1}^{n} (U^{\gamma}y_i \cap A ).$$

Now corresponding to each $y_i$, let $U^{\gamma}xy_i$ be the $\gamma$-open set containing x, then

$U^{\gamma}_x = \cap U^{\gamma}xy_i$ is $\gamma$-open containing x, since $\gamma$ is regular. Also

$$U^{\gamma}_x \cap A = U^{\gamma}_x \cap (\bigcup_{i=1}^{n} (U^{\gamma}y_i \cap A )) \subseteq U^{\gamma}_x \cap (\bigcup_{i=1}^{n} U^{\gamma}y_i )$$

$$\subseteq \bigcup_{i=1}^{n} (U^{\gamma}_x \cap U^{\gamma}y_i )$$

$$= \varnothing$$

or $U^{\gamma}_x \cap A = \varnothing$. Hence $U^{\gamma}_x \subseteq X - A$ implies $x \in int_{\gamma} (X - A)$. Consequently,

$X - A = int_{\gamma} (X - A)$. That is, $X - A$ is $\gamma$-open. So, A is $\gamma$-closed. This completes the proof.

**Definition 3**. Let X be a space and $A \subseteq X$. Then the class of $\gamma$-open sets in A is defined in a natural way as :

$$\tau\gamma_A = \{ A \cap O : O \in \tau_{\gamma} \} ,$$

where $\tau_{\gamma}$ is the class of $\gamma$-open sets of X. That is, G is $\gamma$-open in A iff $G = A \cap O$, where O is a $\gamma$-open set in X.

**Definition** [8]. An operation $\gamma$ on $\tau$ is said to be regular, if for any open nbds U , V of $x \in X$, there exists an open nbd W of x such that $U^{\gamma} \cap V^{\gamma} \supseteq W^{\gamma}$.

## 2. $\gamma^s$-regular space

**Definition 4**. A space X is said to be $\gamma^s$-regular space, if for any closed set A and $x \notin A$, there exist open sets U , V such that $x \in U$, $A \subseteq V$ and $U^{\gamma} \cap V^{\gamma} = \varnothing$.

**Example**. Let X= {a,b,c}, $\tau$ ={ $\varnothing$, X, {a}, {b,c}}.

For b$\in$X, define an operation $\gamma : \tau \longrightarrow$ P(X) by

$$\gamma(A) = \begin{cases} \text{A}, & \text{if b} \in A \\ \text{cl(A)}, & \text{if b} \notin A \end{cases}$$

Then easy calculations show that X is a $\gamma^s$-regular space.

**Theorem 3**. Every subspace of $\gamma^s$-regular space X is $\gamma^s$-regular, where $\gamma$ is regular.

**Proof** . Let Y be a subspace of a $\gamma^s$-regular space X. Suppose A is $\gamma$-closed set in Y and y$\in$Y such that y$\notin$A. Then A = B $\cap$ Y, where B is $\gamma$-closed in X. Then y$\notin$B. Since X is $\gamma^s$-regular, there exist open sets U , V in X such that y$\in$U, B $\subseteq$ V and $U^\gamma \cap V^\gamma = \varnothing$. Then U $\cap$ Y and V $\cap$ Y are open sets in Y containing y and A respectively, also

$$\begin{aligned} (U \cap Y)^\gamma \cap (V \cap Y)^\gamma &\subseteq ( U^\gamma \cap Y^\gamma ) \cap ( V^\gamma \cap Y^\gamma ) \quad (\gamma \text{ is regular}) \\ &= ( U^\gamma \cap V^\gamma ) \cap Y^\gamma \\ &= \varnothing \cap Y^\gamma = \varnothing. \end{aligned}$$

This completes the proof.

## 3. $\gamma^s$-normal spaces

**Definition 5**. A space X is said to be $\gamma^s$ -normal space, if for any disjoint closed sets A , B of X, there exist open sets U, V such that A $\subseteq$ U, B $\subseteq$ V and $U^\gamma \cap V^\gamma = \varnothing$.

**Example** . Let X= {a,b,c,d},

$\tau$ ={ $\varnothing$, X, {a}, {b}, {a,b}, {b,d},{a,b,d},{b,c},{b,c,d},{a,b,c}}.

For b$\in$X, define an operation $\gamma : \tau \longrightarrow$ P(X) by

$$\gamma(\text{A}) = \begin{cases} \text{cl(A)}, & \text{if b} \in A \\ \text{clintcl(A)}, & \text{if b} \in A \end{cases}$$

Then X is $\gamma^s$-normal.

Next, we characterize $\gamma^s$ -normal space as :

**Theorem 4**. A space X is $\gamma^s$ -normal if for any closed set A and open set U containing A, there is an open set V containing A such that

$$A \subseteq V \subseteq cl_\gamma (V^\gamma) \subseteq U^\gamma,$$

where $\gamma$ is $\gamma$-open and strictly regular.

**Proof** . Let A , B be disjoint closed sets in X. Then $A \subseteq X–B$, where X–B is open in X. By hypothesis, there is a open set V such that

$$cl_\gamma (V^\gamma) \subseteq ( X– B)^\gamma \qquad \ldots \qquad (1)$$

(1) gives $B^\gamma \subseteq (X–cl_\gamma (V))^\gamma$ and $V\cap (X– cl_\gamma (V^\gamma)) = \varnothing$. Consequently, $A \subseteq V$, $B \subseteq X– cl_\gamma (V^\gamma)$ and $V^\gamma \cap ((X– cl_\gamma (V^\gamma)))^\gamma = \varnothing$.

This proves that X is $\gamma$-normal. This completes the proof.

**Theorem 5**. A $\gamma^s$ -normal $\gamma$-$T_1$ space is $\gamma^s$ –regular, where $\gamma$ is strictly regular.

**Proof** . Suppose A is a closed set and $x\notin A$. Since X is a $\gamma$-$T_1$ space, therefore by Proposition 4.9 [8], each $\{x\}$ is $\gamma$-closed in X. Since X is $\gamma^s$-normal, therefore there exist open sets U , V such that $\{x\} \subseteq U$, $A \subseteq V$ and $U \cap V = \varnothing$, or $x \in U$, $A \subseteq V$ and $U \cap V = \varnothing$ implies that $U^\gamma \cap V^\gamma = \varnothing$, since $\gamma$ is strictly regular. Thus X is $\gamma^s$ -regular.This completes the proof.

**Theorem 6**. A closed subspace of a $\gamma^s$ -normal space X is $\gamma^s$ -normal, where $\gamma$ is regular.

**Proof** . Let A be a closed subspace of $\gamma^s$ -normal space X. Let $A_1$, $A_2$ be disjoint closed sets of A. Then there are closed sets $B_1$, $B_2$ in X such that $A_1= B_1 \cap A$ , $A_2 = B_2 \cap A$. Since A is closed in X, therefore $A_1$, $A_2$ are closed in X. Since X is $\gamma^s$-normal, there exist open sets $U_1$, $U_2$ in X such that $A_1 \subseteq U_1$, $A_2 \subseteq U_2$ and $U^{\gamma}_1 \cap U^{\gamma}_2 = \varnothing$. But then $A_1 \subseteq A \cap U_1$, $A_2 \subseteq A \cap U_2$ , Where $A \cap U_1$, $A \cap U_2$ are open in A and

$$(A \cap U_1)^{\gamma} \cap (A \cap U_2)^{\gamma} \subseteq (A^{\gamma} \cap U^{\gamma}_1) \cap (A^{\gamma} \cap U^{\gamma}_2) \quad \text{( since } \gamma \text{ is regular)}$$

$$= A^{\gamma} \cap (U^{\gamma}_1 \cap U^{\gamma}_2)$$

$$= A^{\gamma} \cap \varnothing = \varnothing.$$

This proves that A is $\gamma^s$ -normal. Hence the proof.

**Theorem 7**. Every $\gamma_0$-compact and $\gamma$-$T_2$ space is $\gamma^s$-normal, where $\gamma$ is regular and $\gamma$-open.

**Proof**. Let X be $\gamma_0$-compact and $\gamma$-$T_2$ space and $C_1$, $C_2$ be ant two disjoint $\gamma$-closed subsets of X. Then being $\gamma$-closed subset of $\gamma_0$-compact space, $C_1$ is $\gamma_0$-compact. By Theorem 1, for $\gamma_0$-compact $C_2$ and $x \notin C_2$, there are open sets $U_x$, $V_x$ such that

$$x \in U^{\gamma}_x,\ C_2 \subseteq V^{\gamma}_x \text{ and } U^{\gamma}_x \cap V^{\gamma}_x = \varnothing. \qquad \ldots \qquad (1)$$

Let the set $\{ U^{\gamma}_x : x \in C \}$ be a cover of $C_1$ by $\gamma$-open and $\gamma$-closed sets of $C_1$. Since $C_1$ is $\gamma_0$-compact, so there are finite number of elements $x_1, x_2, \ldots, x_n$ such that

$$C_1 \subseteq \bigcup_{i=1}^{n} cl_{\gamma} (U^{\gamma}x_i) = \bigcup_{i=1}^{n} (U^{\gamma}x_i).$$

Let $U = \bigcup_{i=1}^{n} (U^{\gamma}x_i)$ , $V = \bigcap_{i=1}^{n} (V^{\gamma}x_i)$. Then $C_1 \subseteq U$, $C_2 \subseteq V$ and

$$(U \cap V)^{\gamma} = ((\bigcup_{i=1}^{n} (U^{\gamma}x_i) \cap (\bigcap_{i=1}^{n} (V^{\gamma}x_i))^{\gamma}$$

$$= (\varnothing)^{\gamma} = \varnothing.$$

Hence X is $\gamma^s$-normal. This completes the proof.